\documentclass[12pt]{article}

\usepackage{latexsym,amsmath}
\usepackage{amssymb}
\usepackage{url}
\usepackage{graphicx}
\usepackage{microtype}

\newcommand{\cE}[0]{{\mathcal E}}

\newcommand{\cM}[0]{{\mathcal M}}

\newcommand{\cO}[0]{{\mathcal O}}

\newcommand{\cS}[0]{{\mathcal S}}

\newtheorem{thm}{Theorem}[section]

\newtheorem{lemma}[thm]{Lemma}

\newtheorem{cor}[thm]{Corollary}

\newcommand{\gl}[0]{\lambda}

\newcommand{\grg}[0]{\gamma}

\begin{document}

\renewcommand{\thefootnote}{\fnsymbol{footnote}}

\title{$H$-colouring bipartite graphs}

\author{John Engbers
\and
David Galvin}

\date{\today\thanks{$\{$jengbers, dgalvin1$\}$@nd.edu; Department of Mathematics,
University of Notre Dame, Notre Dame IN 46556. Galvin in part supported by National Security Agency grant H98230-10-1-0364.}}

\maketitle

\begin{abstract}

For graphs $G$ and $H$, an {\em $H$-colouring} of $G$ (or {\em homomorphism} from $G$ to $H$) is a function from the vertices of $G$ to the vertices of $H$ that preserves adjacency. $H$-colourings generalize such graph theory notions as proper colourings and independent sets.

For a given $H$, $k \in V(H)$ and $G$ we consider the proportion of vertices of $G$ that get mapped to $k$ in a uniformly chosen $H$-colouring of $G$. Our main result concerns this quantity when $G$ is regular and bipartite. We find numbers $0 \leq a^-(k) \leq a^+(k) \leq 1$ with the property that for all such $G$, with high probability the proportion is between $a^-(k)$ and $a^+(k)$, and we give examples where these extremes are achieved.
For many $H$ we have $a^-(k) = a^+(k)$ for all $k$ and so in these cases we obtain a quite precise description of the almost sure appearance of a randomly chosen $H$-colouring.

As a corollary, we show that in a uniform proper $q$-colouring of a regular bipartite graph, if $q$ is even then with high probability every colour appears on a proportion close to $1/q$ of the vertices, while if $q$ is odd then with high probability every colour appears on at least a proportion close to $1/(q+1)$ of the vertices and at most a proportion close to $1/(q-1)$ of the vertices.

Our results generalize to natural models of weighted $H$-colourings, and also to bipartite graphs which are sufficiently close to regular. As an application of this latter extension we describe the typical structure of $H$-colourings of graphs which are obtained from $n$-regular bipartite graphs by percolation, and we show that $p=1/n$ is a threshold function across which the typical structure changes.

The approach is through entropy, and extends work of J. Kahn, who considered the size of a randomly chosen independent set of a regular bipartite graph.

\end{abstract}

\section{Introduction and statement of results}
\label{sec-intro}

Let $G=(V(G),E(G))$ be a simple, loopless, finite graph, and let $H=(V(H),E(H))$ be a finite graph without multiple edges but perhaps with loops. An {\em $H$-colouring} of $G$, or {\em homomorphism} from $G$ to $H$, is a function from $V(G)$ to $V(H)$ that preserves adjacency. The set of $H$-colourings of $G$ is thus
$$
{\rm Hom}(G,H)=\{f:V(G)\rightarrow V(H): uv \in E(G) \Rightarrow f(u)f(v) \in E(H)\}.
$$
$H$-colourings generalize a number of important graph theory notions. For example, when $H$ is the complete graph on $q$ vertices, ${\rm Hom}(G,H)$ coincides with the set of proper $q$-colourings of $G$, and when $H$ consists of two vertices joined by an edge, with a loop at one of the vertices, then ${\rm Hom}(G,H)$ may be identified with the set of independent sets of $G$, via the preimage of the unlooped vertex.

$H$-colourings have a natural statistical physics interpretation as configurations in {\em hard-constraint spin models}. Here, the vertices of $G$ are thought of as sites that are occupied by particles, with edges of $G$ representing pairs of bonded sites. The vertices of $H$ are the different types of particles (or spins), and the occupation rule is that bonded sites must be occupied by pairs of particles that are adjacent in $H$. A legal configuration in such a spin model is exactly an $H$-colouring of $G$.

From the statistical physics standpoint, there is a very natural family of probability distributions that can be put on ${\rm Hom}(G,H)$. Fix a set of positive weights $\Lambda = \{\gl_i:i \in V(H)\}$ indexed by the vertices of $H$. We think of the magnitude of $\lambda_k$ as measuring how likely particle $k$ is to appear at each site. This can be formalized by giving each $f \in {\rm Hom}(G,H)$ weight $w_\Lambda(f) = \prod_{v \in V(G)} \gl_{f(v)}$ and probability
$$
p_\Lambda(f) = \frac{w_\Lambda(f)}{Z_\Lambda(G,H)}
$$
where $Z_\Lambda(G,H)=\sum_{f \in {\rm Hom}(G,H)} w_\Lambda(f)$ is the appropriate normalizing constant or {\em partition function} of the model. For an introduction to statistical physics spin models from a combinatorial perspective, see for example \cite{BrightwellWinkler2}.

The question to be addressed in this paper is the following. What can be said about an $f$ that is drawn from ${\rm Hom}(G,H)$ according to the distribution $p_\Lambda$? Specifically, for each $f \in {\rm Hom}(G,H)$ and $k \in V(H)$ set
$$
s(k,f) = \frac{|f^{-1}(k)|}{|V(G)|}
$$
and
$$
\bar{p}_\Lambda(k) = \frac{1}{|V(G)|} \sum_{v \in V(G)} p_\Lambda(f(v)=k)~\left(= E_\Lambda\left(s(k,f)\right)\right).
$$
The aim of this paper is to give fairly precise estimates for $\bar{p}_\Lambda(k)$ and the distribution of $s(k,f)$ for $f$ chosen according to $p_\Lambda$, when $G$ is bipartite and either regular or sufficiently close to regular.

\medskip

The point of departure for this work is a result of Kahn on the hard-core model. When $H=H_{\rm ind}$ with $V(H_{\rm ind})=\{0,1\}$ and $E(H_{\rm ind})=\{00,01\}$, the set of vertices of $G$ mapped to $1$ forms an independent set in $G$, and ${\rm Hom}(G,H_{\rm ind})$ can be identified with ${\mathcal I}(G)$, the set of independent sets in $G$. For each $\gl > 0$, the {\em hard-core} model on $G$ is the probability distribution ${\rm hc}(\lambda)$ on ${\mathcal I}(G)$ that assigns to each $I \in {\mathcal I}(G)$ a probability proportional to $\gl^{|I|}$. One of the oldest and most studied spin models in statistical physics, this is a simple mathematical model of the occupation of space (represented by $G$) by particles of non-negligible size. The model can easily be realized as a spin model of the kind described above by assigning weights $\lambda_0=1$ and $\lambda_1=\gl$ to the vertices of $H_{\rm ind}$.

Kahn \cite{Kahn} studied this model on a regular bipartite graph $G$. He proved that for all fixed $\lambda > 0$, the model exhibits a {\em phase coexistence} in the sense that if $G$ has equipartition ${\mathcal E} \cup {\mathcal O}$ then most ${\rm hc}(\gl)$ independent sets tend to come either mostly from ${\mathcal E}$ or mostly from ${\mathcal O}$, in the sense that the size of an independent set chosen according to ${\rm hc}(\gl)$ is concentrated close to $\gl/(2(1+\gl$)), which is exactly the expected size of an independent set chosen according to the distribution that half the time picks a ${\rm hc}(\gl)$ independent set from ${\mathcal E}$ and half the time picks from ${\mathcal O}$. The following theorem (\cite[Theorem 1.4 \& Corollary 1.5]{Kahn}) formalizes this.
\begin{thm} \label{thm-Kahn-indsetsgeneral}
Let $\gl > 0$ be fixed. There are positive constants $c_1$, $c_2$, $c_3$ and $c_4$ (depending on $\gl$) such that for every $d$-regular bipartite graph $G$ on $N$ vertices, the following two statements hold. Firstly, for every $\varepsilon \geq c_1/\sqrt{d}$, if $I$ is chosen from ${\mathcal I}(G)$ according to the distribution ${\rm hc}(\lambda)$ then
$$
{\rm Pr}\left(\left||I|- \frac{\lambda N}{2(1+\lambda)}\right| \geq \varepsilon N\right) \leq c_2\varepsilon^{-1}2^{-c_3\varepsilon^2 N}.
$$
Secondly,
$$
\left|\frac{E(|I|)}{N} - \frac{\gl}{2(1+\gl)}\right| \leq c_4\zeta
$$
where
\begin{equation} \label{def-zeta}
\zeta = \max\left\{\frac{1}{\sqrt{d}},\sqrt{\frac{\log N}{N}}\right\}.
\end{equation}
\end{thm}
In particular, a uniformly chosen independent set ($\gl=1$) from a regular bipartite graph consists, with high probability, of close to one quarter of the vertices. While this corollary may seem more natural than the formulation of Theorem \ref{thm-Kahn-indsetsgeneral}, it is worth noting that in order to prove the theorem in the special case of $\gl=1$ it is necessary (at least using the entropy methods of \cite{Kahn}) to pass to the more general weighted model first. Similarly, it might seem more natural in the present paper to focus on the structure of uniform $H$-colourings, but we are unable to obtain any results without introducing weights.

From (\ref{def-zeta}) we see that Theorem \ref{thm-Kahn-indsetsgeneral} only gives a concentration result when we consider families of graphs with $d$ going to infinity. This is not just an artifact of the proof. For families of graphs with $d$ fixed (and only $N$ going to infinity), the behavior of $E(|I|)/N$ depends very much on the particular choice of family. As an example, consider the case $d=2$. If $G_N$ is the disjoint union of $N/4$ copies of the cycle $C_4$, and $I$ is chosen uniformly from $I$ is chosen uniformly from ${\mathcal I}(G)$, then $E(|I|)/N$ is easily seen to be concentrated close to $2/7$. If, however, $G_N$ is the disjoint union of $N/6$ copies of the cycle $C_6$, then $E(|I|)/N$ concentrated close to $5/18$. For this reason we implicitly assume throughout that $d$ going to infinity.

\medskip

We now set up some notation that allows us to state our main result, which is an extension of Theorem \ref{thm-Kahn-indsetsgeneral} to arbitrary weighted $H$-colourings.
From now on, whenever $H$ and $\Lambda$ are mentioned, it will be assumed that $H$ is a finite graph without multiple edges but perhaps with loops, and that $\Lambda$ is a set of positive weights indexed by the vertices of $H$.
For $A, B \subseteq V(H)$ write $A \sim B$ if for all $u \in A$ and $v \in B$ we have $uv \in E(H)$, and set
$$
\eta_\Lambda(H) = \max\left\{w_\Lambda(A)w_\Lambda(B):A \sim B\right\}
$$
where
$w_\Lambda(\cdot)=\sum_{i \in \cdot} \gl_i$. Then set
$$
{\mathcal M}_\Lambda(H) = \{(A,B) \in V(H)^2: A \sim B,~w_\Lambda(A)w_\Lambda(B) = \eta_\Lambda(H)\}.
$$
Next define
$$
a_\Lambda^+(k) = \frac{\max\left\{w_\Lambda(A)\gl_k \mathbf{1}_{\{k \in B\}} + w_\Lambda(B)\gl_k \mathbf{1}_{\{k \in A\}} : (A,B) \in \cM_\Lambda(H)\right\}}{2 \eta_\Lambda(H)}
$$
and define $a_\Lambda^-(k)$ similarly, with $\max$ replaced by $\min$. Note that if $k$ does not appear in any $(A,B)\in {\mathcal M}_\Lambda(H)$ then $a_\Lambda^+(k)=0$ and that if there is a pair $(A,B)\in {\mathcal M}_\Lambda(H)$ in which $k$ does not appear then $a_\Lambda^-(k)=0$. Note also that $a_\Lambda^-(k) \leq a_\Lambda^+(k)$. Finally, note that $a_\Lambda^+(k)$ and $a_\Lambda^-(k)$ both take the form
$$
\frac{\gl_k \mathbf{1}_{\{k \in A\}}}{2w_\Lambda(A)} + \frac{\gl_k \mathbf{1}_{\{k \in B\}}}{2w_\Lambda(B)}
$$
for some $(A,B) \in \cM_\Lambda(H)$. We may interpret this quantity as the expected proportion of vertices mapped to $k$ in a $p_\Lambda$-chosen $H$-colouring subject to the condition that all vertices from one partition class of $G$ get mapped to $A$ and all from the other class get mapped to $B$; we will refer to such a colouring as a {\em pure}-$(A,B)$ colouring. Finally, for every $\varepsilon > 0$ and $k \in V(H)$ define  
$$  %
I_{k}(\varepsilon) = [0, a_\Lambda^-(k) - \varepsilon) \cup (a_\Lambda^+(k) + \varepsilon, 1].  %
$$  %

\medskip

Before stating our main result, we motivate it by considering weighted $H$-colourings of $K_{d,d}$, the complete bipartite graph with $d$ vertices in each class, for some fixed $H$ and $\Lambda$. The adjacency structure of $K_{d,d}$ ensures that all $H$-colourings are pure-$(A,B)$ for some $(A,B)$ with $A\sim B$, and that moreover all but a vanishing proportion (in $d$) of $Z_\Lambda(K_{d,d},H)$ comes from pure-$(A,B)$ colourings for some $(A,B) \in \cM_\Lambda(H)$. It follows that for each $k \in V(H)$, in an $H$-colouring chosen according to $p_\Lambda$ we have that with probability $1-o(1)$ the proportion of vertices of $K_{d,d}$ mapped to $k$ will be between $a_\Lambda^-(k) - o(1)$ and $a_\Lambda^+(k) + o(1)$. Our main result, which we now state, asserts that this property of $K_{d,d}$ is essentially shared by all $d$-regular graphs.
\begin{thm} \label{thm-unpeaked}
Fix $H$ and $\Lambda$. 
There are positive constants $c_1$, $c_2$, $c_3$ and $c_4$ (depending on $H$ and $\Lambda$) such that for every $d$-regular bipartite graph $G$ on $N$ vertices, the following two statements hold. Firstly, for every $\varepsilon \geq c_1/\sqrt{d}$ and $k \in V(H)$ we have
\begin{equation} \label{main-prob-bound}
p_\Lambda\left(s(k,f) \in I_k(\varepsilon)\right) \leq c_2\varepsilon^{-1}2^{-c_3\varepsilon^2N}.
\end{equation}
Secondly, for each $k \in V(H)$ we have
\begin{equation} \label{cor-prob-bound}
\bar{p}_\Lambda(k)  \in \left[a_\Lambda^-(k) - c_4\zeta, a_\Lambda^+(k) + c_4\zeta\right]
\end{equation}
where $\zeta$ is as defined in (\ref{def-zeta}).
\end{thm}
In other words, for regular bipartite $G$ the distribution $p_\Lambda$ is concentrated on $H$-colourings for which, for every $k \in V(H)$, the proportion of vertices mapped to $k$ is roughly between $a_\Lambda^-(k)$ and $a_\Lambda^+(k)$.

The proof of Theorem \ref{thm-unpeaked} goes along the following lines. We upper bound the contribution to $Z_\Lambda(G,H)$ from those $f \in {\rm Hom}(G,H)$ with $|f^{-1}(k)|/N = \grg \geq a^+(k)+\varepsilon$ by $Z_{\Lambda(k,\delta)}(G,H)/(1+\delta)^{\grg N}$ for some suitably small $\delta > 0$ (where $\Lambda(k, \delta)$ is obtained from $\Lambda$ by multiplying $\gl_k$ by $1+\delta$ and leaving all other $\gl_i$ unchanged). We in turn upper bound $Z_{\Lambda(k,\delta)}(G,H)$ using a result of Galvin and Tetali \cite{GalvinTetali-weighted} to the effect that for all $H$ and $\Lambda$ and all $d$-regular graphs $G$ on $N$ vertices we have
\begin{equation} \label{from-GT}
Z_\Lambda(G,H) \leq Z_\Lambda(K_{d,d}, H)^\frac{N}{2d}
\end{equation}
where $K_{d,d}$ is the complete bipartite graph with $d$ vertices in each partition class. We upper bound $Z_{\Lambda(k,\delta)}(K_{d,d}, H)$ in terms of $\eta_{\Lambda(k, \delta)}(H)$, and in the end we get, using our choice of $a_\Lambda^+(k)$ and for some sufficiently small $\delta$, an upper bound on the contribution that is significantly smaller than a trivial lower bound on $Z_\Lambda(G,H)$, showing those $f \in {\rm Hom}(G,H)$ with $|f^{-1}(k)|/N \geq a^+(k)+\varepsilon$ do not contribute greatly to the partition function. The same strategy works for $|f^{-1}(k)|/N$ falling significantly below $a^-(k)$. The details (in the more general setting of Theorem \ref{thm-non-regular}) are given in Section \ref{sec-proof-of-main}.

\medskip

When $a_\Lambda^-(k)=a_\Lambda^+(k)$ for all $k$, we obtain a single vector around which $(s(k,f): k \in V(H))$ is concentrated for $f$ chosen according to $p_\Lambda$.
\begin{cor} \label{cor-peaked}
Fix $H$ and $\Lambda$. 
Suppose that for all $k \in V(H)$ there is an $a_\Lambda(k)$ such that $a_\Lambda^-(k)=a_\Lambda^+(k)=a_\Lambda(k)$. Then there are positive constants $c_1$, $c_2$, $c_3$ and $c_4$ (depending on $H$ and $\Lambda$) such that for every $d$-regular, bipartite graph $G$ on $N$ vertices the following two statements hold. Firstly, for $\varepsilon \geq c_1/\sqrt{d}$ we have
$$
p_\Lambda\left(\left|\left|(s(k,f))_{k \in V(H)}-(a_\Lambda(k))_{k \in V(H)}\right|\right|_\infty \geq \varepsilon \right) \leq c_2\varepsilon^{-1}2^{-c_3\varepsilon^2N}.
$$
Secondly, we have
$$
\left|\left|(\bar{p}_\Lambda(k))_{k \in V(H)}-(a_\Lambda(k))_{k \in V(H)}\right|\right|_\infty \leq c_4\zeta
$$
with $\zeta$ as in (\ref{def-zeta}).
\end{cor}
A situation in which Corollary \ref{cor-peaked} applies is when either $\cM_\Lambda(H)=\{(A,A)\}$ or $\cM_\Lambda(H)=\{(A,B),(B,A)\}$ (for some $A \neq B$). This is in a sense the generic situation. Indeed, for every $H$, if the weights $\gl_i$ are chosen from any continuous distribution supported on $\{x \in {\mathbb R}^{|V(H)|}: x > 0\}$, then with probability $1$ we will have $\cM_\Lambda(H)$ of the form described. As we will see in Example C below, Corollary \ref{cor-peaked} also applies in some other natural situations.

\medskip

The gap between $a_\Lambda^-(k)$ and $a_\Lambda^+(k)$ (if there is one) cannot be closed in general, as the first part of the following theorem shows.
\begin{thm} \label{thm-tightbounds}
Fix $H$ and $\Lambda$. 
There is a family $\{G_d\}_{d=1}^\infty$ of $d$-regular bipartite graphs, a function $g(d)=o(1)$ and a positive constant $c$ (depending on $H$ and $\Lambda$) such that for each $k \in V(H)$,
$$
\left.
\begin{array}{r}
p_\Lambda\left( \left|s(k,f)-a_\Lambda^+(k)\right| \leq g(d) \right) \\
p_\Lambda\left( \left|s(k,f)-a_\Lambda^-(k)\right| \leq g(d) \right)
\end{array}
\right\}
  \geq c-g(d).
$$
There is also a family  $\{G'_d\}_{d=1}^\infty$ of $d$-regular bipartite graphs, a function $g(d) = o(1)$ and (for each $k \in V(H)$) an $a_\Lambda(k)$ satisfying $a^-_\Lambda(k) \leq a_\Lambda(k) \leq a^+_\Lambda(k)$ such that for each $k$,
$$
p_\Lambda\left( \left|s(k,f)-a_\Lambda(k)\right| \leq g(d) \right)  \geq 1-g(d)
$$
and
$$
|\bar{p}_\Lambda(k) - a_\Lambda(k)| \leq g(d).
$$
\end{thm}
We prove Theorem \ref{thm-tightbounds} in Section \ref{sec-colourings-proof}. The graphs $G_d$ we exhibit will be suitably chosen random regular graphs, and we will use the expansion of these graphs to show that all but $o(1)$ of $p_\Lambda$ is concentrated on pure-$(A,B)$ colourings for $(A,B) \in \cM_\Lambda(H)$. The graphs $G'_d$ will be disjoint unions of complete bipartite graphs on $2d$ vertices. Basic concentration estimates together with the independence of the components will give the claimed result.

\medskip

We now explore the consequences of Theorem \ref{thm-unpeaked} for some specific choices of $H$ and $\Lambda$.

\medskip

\noindent {\bf Example A} (Hard-core model) Let $H=H_{\rm ind}$ be as described earlier, with $\gl_0=1$ and $\gl_1=\gl$. We have seen that an element of ${\rm Hom}(G,H_{\rm ind})$ chosen according to $p_\Lambda$ is a configuration in the hard-core model on $G$ with activity $\gl$. With these choices we have $\cM_\Lambda(H_{\rm ind})=\{(\{0\},\{0,1\}), (\{0,1\},\{0\})\}$ and
$$
a_\Lambda^-(1)=a_\Lambda^+(1)= \frac{\gl}{2(1+\gl)}
$$
and so Theorem \ref{thm-unpeaked} indeed generalizes Theorem \ref{thm-Kahn-indsetsgeneral}, as claimed.

\medskip

\noindent {\bf Example B} (Multistate hard-core model) Let $H=H_k$ be the graph on vertex set $\{0, \ldots, k\}$ with $ij \in E(H)$ if and only if $i+j \leq k$, and $\gl_i = \gl^i$ for some fixed $\gl>0$. An element of ${\rm Hom}(G,H_k)$ chosen according to $p_\Lambda$ is exactly a configuration of the {\em multistate hard-core} (or {\em multicast communications}) model on $G$ with activity $\gl$. This model allows multiple particles (up to and including $k$) at each site, with the restriction that there are no more than $k$ particles in total across each edge. A generalization of the hard-core model (the case $k=1$), it has been studied in a variety of contexts: in communications \cite{MitraRamananSengupta}, statistical physics \cite{MazelSuhov} and combinatorics \cite{GalvinMartinelliRamananTetali}. For $k$ even the unique pair $(A,B) \in \cM_\Lambda(H_k)$ has $A=B=\{1,\ldots, k/2\}$, while for $k$ odd, say $k=2\ell+1$, we have $\cM_\Lambda(H_k)=\{(A,B),(B,A)\}$ with  $A=\{1, \ldots, \ell\}$ and $B=\{1, \ldots, \ell+1\}$. In either case Corollary \ref{cor-peaked} shows that for this model $(s(k,f): k \in V(H))$  is concentrated close to a single value for $f$ chosen according to $p_\Lambda$.

\medskip

\noindent {\bf Example C} (Uniform proper $q$-colourings) Let $H=K_q$, the complete graph on $q$ vertices, and $\Lambda = (1, \ldots, 1)$. An element of ${\rm Hom}(G,K_q)$ chosen according to $p_\Lambda$ corresponds to a uniform proper $q$-colouring of $G$. In this case elements of $\cM_\Lambda(K_q)$ consist of all partitions of $V(K_q)$ into two classes as near equal in size as possible, and an easy calculation gives that for all colours $k$
$$
a_\Lambda^-(k) = \frac{1}{2 \lceil q/2 \rceil} ~~~\mbox{and} ~~~ a_\Lambda^+(k) = \frac{1}{2 \lfloor q/2 \rfloor}
$$
so that in particular $a_\Lambda^-(k)=a_\Lambda^+(k)= 1/q$ for $q$ even, and we get the following corollary of Theorem \ref{thm-unpeaked}.
\begin{cor} \label{cor-colourings}
Fix $q \in {\mathbb N}$. There are positive constants $c_1$, $c_2$, $c_3$ and $c_4$ (depending on $q$) such that for every $d$-regular, bipartite graph $G$ on $N$ vertices, the following statements hold. If $\chi$ is a uniformly chosen $q$-colouring of $G$ and $\varepsilon \geq c_1/\sqrt{d}$ then for $q$ even
$$
{\rm Pr}\left(\exists k \in V(H): \left|\frac{|\chi^{-1}(k)|}{N} - \frac{1}{q}\right| \geq \varepsilon \right) \leq c_2\varepsilon^{-1}2^{-c_3\varepsilon^2N},
$$
and for $q$ odd
$$
\left.
\begin{array}{r}
{\rm Pr}\left(\exists k \in V(H): \frac{|\chi^{-1}(k)|}{N} \leq \frac{1}{q+1} - \varepsilon\right) \\
{\rm Pr}\left(\exists k \in V(H): \frac{|\chi^{-1}(k)|}{N} \geq \frac{1}{q-1} + \varepsilon\right)
\end{array}
\right\} < c_2\varepsilon^{-1}2^{-c_3\varepsilon^2N}.
$$
\end{cor}
So for even $q$, almost all proper $q$-colourings of a regular bipartite graph are ``almost equitable''. Of course, by the symmetry of $K_q$ we have $E(|\chi^{-1}(k)|)=N/q$ for all $k$ in this case.

\medskip

The condition that $G$ be regular can be relaxed quite a bit; we simply require that $G$ has not too many low degree vertices, that the sum of the degrees of high degree vertices is not too large, and that the difference between the sizes of the partition classes is not too great.
\begin{thm} \label{thm-non-regular}
Fix $H$ and $\Lambda$. 
There are positive constants $c_1$, $c_2$, $c_3$ and $c_4$ (depending on $H$ and $\Lambda$) such that the following statements hold. Let $G$ be a bipartite graph on $N$ vertices with bipartition classes $\cE$ and $\cO$ (with $|\cO|\geq |\cE|$). Let $d$ be an arbitrary positive parameter. Let $\varepsilon$ satisfy $\varepsilon \geq c_1\sqrt{h(G,d)}$ where
$$
h(G,d) = \frac{1}{d} + \frac{|\{v \in \cE: d(v) < d\}|}{N}  + \frac{|\cO|-|\cE|}{N} + \frac{1}{dN}\sum_{v \in \cO} (d(v)-d)\mathbf{1}_{\{d(v)\geq d\}}.
$$
Then for each $k \in V(H)$ we have (\ref{main-prob-bound}), as well as (\ref{cor-prob-bound}) with now
$$
\zeta = \max \left\{\sqrt{h(G,d)}, \sqrt{\frac{\log N}{N}}\right\}.
$$
\end{thm}
If $G$ is $d$-regular then $h(G,d) = 1/d$ and so Theorem \ref{thm-non-regular} is a generalization of Theorem \ref{thm-unpeaked}. The proof of Theorem \ref{thm-non-regular} follows the same lines as already described for Theorem \ref{thm-unpeaked}, except that we now require a new upper bound on $Z_\Lambda(G,H)$. In Section \ref{sec-non-reg-proof} we modify the entropy-based proof of (\ref{from-GT}) to obtain the following, which is just what we need for Theorem \ref{thm-non-regular}, the proof of which is then given in Section \ref{sec-proof-of-main}. Here $d(v)=|\{u \in V(G): uv \in E(G)\}|$ is the degree of vertex $v$, and we write $w_\Lambda(H)$ for $w_\Lambda(V(H))$.
\begin{thm} \label{thm-improving-GT}
Fix $H$ and $\Lambda$, 
and suppose that $\gl_i>1$ for all $i \in V(H)$. Let $G$ be any bipartite graph on bipartition classes $\cE$ and $\cO$, with $|\cO|\geq |\cE|$, and let $d$ be an arbitrary positive parameter. Then
$$
Z_\Lambda(G,H) \leq w_\Lambda(H)^{|\{w \in \cE:d(w)< d\}|} \prod_{v \in \cO} Z_\Lambda(K_{d(v),d},H)^\frac{1}{d}.
$$
\end{thm}
Note that if $G$ is $d$-regular then Theorem \ref{thm-improving-GT} reduces to (\ref{from-GT}). Note also that the condition imposed on the $\gl_i$ by Theorem \ref{thm-improving-GT}  is not restrictive: if $\Lambda'$ is obtained from $\Lambda$ by multiplying all $\gl_i \in \Lambda$ by the same positive constant then $p_\Lambda(N_1(f)=\cdot) = p_{\Lambda'}(N_1(f)=\cdot)$ and so we may assume without loss of generality that $\min\{\gl_i : i \in V(H)\}>1$.

\medskip

Theorem \ref{thm-non-regular} is only of interest in situations where $h(G,d)$ can be shown to be small (as, for example, when $G$ is $d$-regular).
A natural situation where we can say something about $h(G,d)$ is in percolation. Given a graph $G$ and a parameter $0 \leq p \leq 1$, let $G_p$ be a random subgraph of $G$ obtained by deleting each edge independently with probability $1-p$ (so the probability that $G_p=H$ is $p^{|E(H)|}(1-p)^{|E(G)|-|E(H)|}$). A corollary of Theorem \ref{thm-non-regular} (which we will prove in Section \ref{sec-perc-proof}) is the following ``phase transition'' phenomenon for percolation on a regular bipartite graph. If $G$ is an $n$-regular bipartite graph and $p$ is much greater than $1/n$, then the typical appearance of a $p_\Lambda$-chosen $H$-colouring of $G_p$ is similar to that of a $p_\Lambda$-chosen $H$-colouring of $G$, whereas if $p$ is much smaller than $1/n$, then as long as there is some $k \in V(H)$ with $\gl_k/w_\Lambda(H) \not \in [a_\Lambda^-(k), a_\Lambda^+(k)]$, these two objects have different appearances.
\begin{cor} \label{cor-percolation}
Fix $H$ and $\Lambda$. 
Let $f(n) =\omega(1)$. There is a function $g(n) = o(1)$ (depending on $f(n)$) such that if $\{G^n\}_{n=1}^{\infty}$ is a sequence of $n$-regular bipartite graphs and $p$ satisfies $p \geq f(n)/n$, then with probability at least $1-g(n)$ the graph $G^n_p$ satisfies that for each $k \in V(H)$ we have
$$
p_\Lambda\left(s(k,f) \in I_k(g(n))\right) \leq g(n)
$$
and
$$
\bar{p}_\Lambda(k) \in \left[a^-_\Lambda(k)-g(n), a^+_\Lambda(k)+g(n)\right].
$$
If on the other hand $p \leq 1/(f(n)n)$ then with probability at least $1-g(n)$ we have that for each $k \in V(H)$,
$$
p_\Lambda\left(\left|s(k,f) - \frac{\gl_k}{w_\Lambda(H)}\right| \leq g(n)\right) \geq 1 - g(n).
$$
and
$$
\left|\bar{p}_\Lambda(k) - \frac{\gl_k}{w_\Lambda(H)}\right| \leq g(n).
$$
\end{cor}
For the multicast model (Example B), for example, we have
$$
a_\Lambda^-(0)=a_\Lambda^+(0) = \frac{1}{2\left(\sum_{i \leq \lfloor k/2 \rfloor} \gl^i\right)} + \frac{1}{2\left(\sum_{i \leq \lceil k/2 \rceil} \gl^i\right)} > \frac{1}{\sum_{i \leq k} \gl^i}
$$
and so Corollary \ref{cor-percolation} shows a phase transition for this model. For the uniform $q$-colouring model (Example C), on the other hand, Corollary \ref{cor-percolation} gives no information about what happens as $p$ crosses $1/n$.

\section{Proof of Theorem \ref{thm-improving-GT}} \label{sec-non-reg-proof}

We will initially assume that for all $i \in V(H)$, we have $\gl_i \in {\mathbb Q}$.
Under this assumption, we can relate $Z_\Lambda(G,H)$ to a uniform model. We repeat an idea used in \cite{GalvinTetali-weighted} and first introduced in \cite{BrightwellWinkler}. Let $C$ be any positive integer with the property that $C\gl_i \in {\mathbb Z}$ for each $i \in V(H)$. Let $H_\Lambda^C$ be the graph obtained from $H$ by the following process: replace each vertex $i$ with a set $S_i$ of size $C\gl_i$, replace each edge $ij$ ($i \neq j$) with a complete bipartite graph between $S_i$ and $S_j$, and replace each loop $ii$ with a complete looped graph on $S_i$. It is easy to check that for any $N$ vertex graph $G$ we have
\begin{equation} \label{ent0}
Z_\Lambda(G,H) = \frac{|{\rm Hom}(G,H_\Lambda^C)|}{C^N}.
\end{equation}
We now bound $|{\rm Hom}(G,H_\Lambda^C)|$ using an entropy approach that was used in \cite{Kahn} to upper bound the number of independent sets in a regular bipartite graph, and was generalized in \cite{GalvinTetali-weighted} to bound $|{\rm Hom}(G,H)|$ for arbitrary $H$ and regular bipartite $G$. We very briefly review the necessary entropy background here; see for example \cite{McEliece} for a more detailed treatment.

For a discrete 
random variable $X$, let $R(X)$ be the support of the mass function of $X$.  
Define the {\em entropy} of $X$ to be
$$
H(X) = \sum_{x \in R(X)} -P(X=x) \log P(X=x),
$$
where here, and throughout the rest of this paper, logarithms have base 2. 
We may think of $H(X)$ as a measure of the randomness of $X$ or as the amount of information it contains. The {\em conditional entropy} of $X$ given the discrete random variable $Y$ is given by
$$
H(X|Y) = \sum_{y \in R(Y)} P(Y=y) \sum_{x \in R(X)} - P(X=x|Y=y) \log P(X=x|Y=y).
$$
Here are the basic facts about the entropy function that we will need. The inequality that makes entropy useful as a tool for enumeration is
\begin{equation} \label{ent-range}
H(X) \leq \log |R(X)|
\end{equation}
with equality if and only if $X$ is uniform.
For a vector $(X_1, \ldots, X_n)$ of random variables (itself a discrete and finite valued random variable) we have a chain rule
\begin{equation} \label{ent-chain}
H(X_1, \ldots, X_n) = H(X_1) + H(X_2 | X_1)  + \ldots + H(X_n | X_1, \ldots, X_{n-1}).
\end{equation}
For random variables $X$, $Y$ and $Z$ we have
\begin{equation} \label{ent-conditioning}
H(X|Y) \leq H(X)~~~\mbox{and}~~~ H(X|Y,Z) \leq H(X|Y)
\end{equation}
(so dropping conditioning does not decrease entropy). Finally, we have conditional subadditivity:
\begin{equation} \label{ent-condsub}
H(X_1, \ldots, X_n|Y) \leq  H(X_1|Y) + H(X_2|Y)  + \ldots + H(X_n|Y).
\end{equation}

Now let $f$ be a uniformly chosen element of ${\rm Hom}(G,H_\Lambda^C)$. By (\ref{ent-chain}) the entropy of $f$ satisfies
\begin{equation} \label{ent1}
H(f) = H(f(\cE)) + H(f(\cO)|f(\cE)).
\end{equation}
We upper bound $H(f(\cO)|f(\cE))$ using (\ref{ent-conditioning}) and (\ref{ent-condsub}):
\begin{equation} \label{ent2}
H(f(\cO)|f(\cE)) \leq \sum_{v \in \cO} H(f(v)|f(N(v)))
\end{equation}
where $N(v)=\{u \in V(G): uv \in E(G)\}$ is the neighbourhood of $v$. We upper bound $H(f(\cE))$ using a form of Shearer's Lemma \cite{ChungFranklGrahamShearer} derived from Radhakrishnan's proof of same (see for example \cite{Kahn2}). Put a total order $<$ on the vertices of $G$. For each $v \in \cO$ with $N(v)=\{n_1, \ldots, n_{d(v)}\}$ where $n_1 < \ldots < n_{d(v)}$ we have, by (\ref{ent-chain}) and (\ref{ent-conditioning}),
\begin{eqnarray}
H(f(N(v))) & = & \sum_{i=1}^{d(v)} H(f(n_i)|f(n_{i-1}), \ldots, f(n_1)) \nonumber \\
& \geq & \sum_{i=1}^{d(v)} H(f(n_i)|\{f(u): u < n_i\}) \nonumber
\end{eqnarray}
and so
\begin{eqnarray}
\sum_{v \in \cO} H(f(N(v)))  & \geq  & \sum_{w \in \cE} d(w)H(f(w)|\{f(u):u < w\}) \nonumber \\
& = & \sum_{w \in \cE} (d + (d(w)-d))H(f(w)|\{f(u):u < w\}) \label{ent3}
\end{eqnarray}
where $d$ is any positive parameter. Since by (\ref{ent-chain}) again we have
$$
\sum_{w \in \cE} H(f(w)|\{f(u):u < w\}) = H(f(\cE))
$$
we rearrange the terms of (\ref{ent3}) to get
\begin{equation} \label{ent4}
H(f(\cE)) \leq \frac{1}{d} \sum_{v \in \cO} H(f(N(v))) + \sum_{w \in \cE} \left(1-\frac{d(w)}{d}\right)H(f(w)|\{f(u):u < w\}).
\end{equation}
We  combine (\ref{ent1}), (\ref{ent2}) and (\ref{ent4}) to upper bound $H(f)$ as the sum of
\begin{equation} \label{ent5}
\frac{1}{d} \sum_{v \in \cO} \left(H(f(N(v))) + dH(f(v)|f(N(v)))\right)
\end{equation}
and
\begin{equation} \label{ent6}
\sum_{w \in \cE} \left(1-\frac{d(w)}{d}\right)H(f(w)|\{f(u):u < w\}).
\end{equation}
We deal first with (\ref{ent5}).
Fix $v \in \cO$. For each $A \in V(H)^{N(v)}$ that occurs as a value of $f(N(v))$, let $p(A)$ be the probability that $A$ occurs and let $e(A)$ be the number of possible ways of assigning an image to $v$ given that $f(N(v))$ takes value $A$. Expanding out the entropy terms we have
\begin{eqnarray}
H(f(N(v))) + dH(f(v)|f(N(v)) & \leq & \sum_A p(A)\log \frac{e(A)^d}{p(A)} \label{ent7} \\
& \leq & \log \sum_A e(A)^d \label{ent8} \\
& \leq & \log |{\rm Hom}(K_{d(v),d},H_\Lambda^C)| \label{kdd} \\
& = & \log \left( C^{d(v)+d} Z_\Lambda(K_{d(v),d},H)\right). \label{ent21}
\end{eqnarray}
We use (\ref{ent-range}) to obtain (\ref{ent7}) and Jensen's inequality for (\ref{ent8}), and the equality in  (\ref{ent21}) follows from (\ref{ent0}). To see (\ref{kdd}) note that we specify an element of ${\rm Hom}(K_{d(v),d},H_\Lambda^C)$ by first choosing the restriction $A$ of the homomorphism to the partition class of size $d(v)$ and then for each of the remaining $d$ vertices choosing the value independently from $e(A)$.  Summing over $v \in \cO$ we see that (\ref{ent5}) is bounded above by
\begin{equation} \label{ent9}
 \log \left(C^{\frac{|E(G)|}{d} + |\cO|}\right) + \sum_{v \in \cO} \log \left(Z_\Lambda(K_{d(v),d},H)^\frac{1}{d}\right).
\end{equation}
For (\ref{ent6}), if $d(w) < d$, we upper bound 
$$
H(f(w)|\{f(u):u < w\}) \leq \log |V(H^C_\Lambda)| = \log \left(Cw_\Lambda(H)\right)
$$
using (\ref{ent-range}) and (\ref{ent-conditioning}).
If $d(w) \geq d$ then we need a lower bound on $H(f(w)|\{f(u):u < w\})$. Since $f$ is a homomorphism, there is at least one $i$ such that $f$ can take values in $S_i$. Fix one such. If we add the condition that $f(w)\in S_i$ then, because the vertices of $S_i$ are indistinguishable, $f(w)$ becomes uniform and its entropy is the logarithm of $|S_i|$. That is,
\begin{eqnarray*}
H(f(w)|\{f(u):u < w\}) & \geq & H(f(w)|\{f(u):u < w\}, \{f(w) \in S_i\}) \nonumber \\
& = & \log C\gl_i \nonumber \\
& \geq & \log C, \nonumber
\end{eqnarray*}
the last inequality following from $\gl_i>1$ for all $i \in V(H)$.
It follows that (\ref{ent6}) is bounded above by
\begin{equation} \label{ent11}
\log \left(C^{|\cE| - \frac{|E(G)|}{d}} w_\Lambda(H)^{|\{w \in \cE:d(w) < d\}|}\right).
\end{equation}
Putting (\ref{ent9}) and (\ref{ent11}) into (\ref{ent1}), using $H(f)=\log |{\rm Hom}(G, H_\Lambda^C)|$ (since $f$ is uniform) and combining with (\ref{ent0}), we obtain Theorem \ref{thm-improving-GT} for rational $\gl_i$'s.
By continuity, this bound remains valid when the $\gl_i$'s are not necessarily rational.

\section{Proof of Theorem \ref{thm-non-regular}}
\label{sec-proof-of-main}

We begin by using Theorem \ref{thm-improving-GT} to put an upper bound on $Z_\Lambda(G,H)$. We first consider those $v \in \cO$ with $d(v) \geq d$. For each of the $4^{|V(H)|}$ ordered pairs $A \sim B$ of subsets of $H$, the contribution to $Z_\Lambda(K_{d(v),d},H)$ from those $f$ with the partition class of $K_{d(v),d}$ of size $d(v)$ mapped to $A$ and the class of size $d$ mapped to $B$ is at most
$$
w_\Lambda(A)^{d(v)}w_\Lambda(B)^d \leq \left(w_\Lambda(A)w_\Lambda(B)\right)^d w_\Lambda(H)^{d(v)-d}
$$
and so
$$
Z_\Lambda(K_{d(v),d},H) \leq 4^{|V(H)|}w_\Lambda(H)^{d(v)-d} \eta_\Lambda(H)^d.
$$
Similarly, for those $v \in \cO$ with $d(v) < d$ we have
$$
Z_\Lambda(K_{d(v),d},H) \leq Z_\Lambda(K_{d,d},H) \leq 4^{|V(H)|}\eta_\Lambda(H)^d.
$$
It follows from Theorem \ref{thm-improving-GT} that $Z_\Lambda(G,H)$ is upper bounded by
$$
\eta_\Lambda(H)^{|\cO|} 4^{\frac{|V(H)||\cO|}{d}} w_\Lambda(H)^{|\{v \in \cE: d(v) < d\}| + \frac{1}{d}\sum_{v \in \cO} (d(v)-d)\mathbf{1}_{\{d(v)\geq d\}}}
$$
and so, using $|\cO| = N/2 + (|\cO|-|\cE|)/2$,
\begin{equation} \label{ub}
Z_\Lambda(G,H) \leq \eta_\Lambda(H)^{\frac{N}{2}} C^{Nh(G,d)}
\end{equation}
where $C$ is a positive constant depending only on $H$ and $\Lambda$.
On the other hand, we get a lower bound (with any $w_\Lambda(A)w_\Lambda(B)=\eta_\Lambda(H)$, and using $\gl_i>1$ for all $i \in V(H)$) by
\begin{eqnarray}
Z_\Lambda(G,H) & \geq & w_\Lambda(A)^{|\cE|}w_\Lambda(B)^{|\cO|} \nonumber \\
& \geq & \left(w_\Lambda(A)w_\Lambda(B)\right)^{|\cE|} \label{lb_int} \\
& = & \eta_\Lambda(H)^{\frac{N}{2}}\eta_\Lambda(H)^{-\frac{|\cO|-|\cE|}{2}}. \label{lb}
\end{eqnarray}
In (\ref{lb_int}) we are using $|\cO|\geq |\cE|$.

We now use (\ref{ub}) and (\ref{lb}) to prove (\ref{main-prob-bound}).
Fix $k \in V(H)$ and an integer $N_k$ satisfying $0  \leq N_k \leq N$ and
$$
\frac{N_k}{N}  \in  [0, a_\Lambda^-(k) - \varepsilon) \cup (a_\Lambda^+(k) + \varepsilon, 1]~\left(= I_k(\varepsilon)\right).
$$
Write $c_k(N_k)$ for the contribution to $Z_\Lambda(G,H)$ from those $f \in {\rm Hom}(G,H)$ with $|f^{-1}(k)|=N_k$. We aim to obtain an upper bound on $c_k(N_k)$ (via (\ref{ub})) which is substantially lower than the lower bound (\ref{lb}), indicating that this term does not contribute greatly to $Z_\Lambda(G,H)$.

We begin by considering $N_k$ for which
$$
\gamma := \frac{N_k}{N}  = a_\Lambda^+(k) + \varepsilon'
$$
for some $\varepsilon'$ satisfying  $\varepsilon \leq \varepsilon' \leq 1 - a_\Lambda^+(k)$.
For any $\delta > 0$ let $\Lambda(k,\delta)$ be obtained from $\Lambda$ by replacing $\gl_k$ with $(1+\delta)\gl_k$ and leaving all other $\gl_i$'s unchanged. By (\ref{ub}) we have
\begin{eqnarray}
(1+\delta)^{N_k}c_k(N_k) & \leq & Z_{\Lambda(k,\delta)}(G,H) \nonumber \\
& \leq & \eta_{\Lambda(k,\delta)}(H)^\frac{N}{2} C^{Nh(G,d)} \label{plus}
\end{eqnarray}
where now the constant $C$ depends on $\delta$ as well as on $H$ and $\Lambda$.

Before proceeding, we need to understand $\eta_{\Lambda(k,\delta)}(H)$. Viewed as a function of $\delta$, the quantity $w_{\Lambda(k,\delta)}(A)w_{\Lambda(k,\delta)}(B)$ (for $(A,B) \in \cM_\Lambda(H)$) is of the form $a+b\delta+c\delta^2$ where $a=\eta_\Lambda(H)$, $b=w_\Lambda(A)\gl_k \mathbf{1}_{\{k \in B\}} + w_\Lambda(B)\gl_k \mathbf{1}_{\{k \in A\}}$ and $c=\gl_k^2 \mathbf{1}_{\{k \in A\cap B\}}$. From this formulation we can easily identify that set $\emptyset \neq \cS^+_\Lambda(k,H) \subseteq \cM_\Lambda(H)$ with the property that for all $\delta > 0$, all $(A,B) \in \cM_\Lambda(H)$ and all $(A',B') \in \cS^+_\Lambda(k,H)$ we have $w_{\Lambda(k,\delta)}(A')w_{\Lambda(k,\delta)}(B') \geq w_{\Lambda(k,\delta)}(A)w_{\Lambda(k,\delta)}(B)$: $\cS^+_\Lambda(k,H)$ consists of all those $(A', B') \in \cM_\Lambda(H)$ for which $b$ is maximum and (subject to this condition) $c$ is maximum. This latter condition simply means that if some of the pairs that maximize $b$ have $c>0$ we only take those pairs, and if they all have $c=0$ we take all pairs.

It is easily seen that there is $\delta^+_k > 0$ (depending on $H$ and $\Lambda$) with the property that for all $0 < \delta < \delta^+_k$ and $(A',B') \in \cS^+_\Lambda(k,H)$ we have $(A',B') \in \cM_{\Lambda(k,\delta)}(H)$. Choose one such, $(A^+,B^+)$, arbitrarily. Note that by construction
$$
a_\Lambda^+(k) = \frac{w_\Lambda(A^+)\gl_k \mathbf{1}_{\{k \in B^+\}} + w_\Lambda(B^+)\gl_k \mathbf{1}_{\{k \in A^+\}}}{2 \eta_\Lambda(H)} = \frac{\gl_k \mathbf{1}_{\{k \in A^+\}}}{2w_\Lambda(A^+)} + \frac{\gl_k \mathbf{1}_{\{k \in B^+\}}}{2w_\Lambda(B^+)}.
$$
Now combining (\ref{lb}) and (\ref{plus}) and choosing $\delta < \delta^+_k$ we have
\begin{eqnarray}
p_\Lambda(|f^{-1}(k)|=N_k) & = & \frac{c_k(N_k)}{Z_\Lambda(G,H)} \nonumber \\
& \leq &  C^{h(G,d)N} \left(\frac{w_{\Lambda(k,\delta)}(A^+)w_{\Lambda(k,\delta)}(B^+)}{w_\Lambda(A^+)w_\Lambda(B^+)(1+\delta)^{2(a_\Lambda^+(k)+\varepsilon')}}\right)^\frac{N}{2} \label{int1}
\end{eqnarray}
where, by our restriction on $\delta$, $C$ may be taken to depend only on $H$ and $\Lambda$.
Our aim is to show that there is a positive constant $c$ (depending on $H$ and $\Lambda$) such that for all $0 < \varepsilon' \leq 1 - a_\Lambda^+(k)$ we can find a $0 < \delta < \delta^+_k$ for which
\begin{equation} \label{to-show}
\frac{w_{\Lambda(k,\delta)}(A^+)w_{\Lambda(k,\delta)}(B^+)}{w_\Lambda(A^+)w_\Lambda(B^+)(1+\delta)^{2(a_\Lambda^+(k)+\varepsilon')}} \leq 2^{-c\varepsilon'^2}.
\end{equation}
Combining this with (\ref{int1}) we see that if $\varepsilon > c\sqrt{h(G,d)}$ for some suitably large positive constant $c$ (depending on $\Lambda$ and $H$) then for all $\varepsilon < \varepsilon' \leq 1 - a_\Lambda^+(k)$ for which $a^+(k)N + \varepsilon' N$ is an integer we have
$$
p_\Lambda\left(|f^{-1}(k)|=a^+(k)N + \varepsilon' N\right)  \leq 2^{-c'{\varepsilon'}^2 N}
$$
for a suitable positive $c'$, and so
\begin{eqnarray}
p_\Lambda\left(|f^{-1}(k)| \geq a^+(k)N + \varepsilon N\right) & \leq & \sum_{\ell \geq \varepsilon N} 2^{-\frac{c'\ell^2}{N}} \nonumber \\
& \leq & 2^{-c'\varepsilon^2 N} \sum_{\ell \geq 0} 2^{-2\ell c'\varepsilon}  \nonumber  \\
& \leq & c''\varepsilon^{-1}2^{-c'\varepsilon^2 N} \label{one}
\end{eqnarray}
for suitably large $c''$ (depending on $c'$). An almost identical argument (the details of which we leave to the reader) yields
\begin{equation} \label{two}
p_\Lambda\left(|f^{-1}(k)| \leq a^-(k)N - \varepsilon N\right) \leq c''\varepsilon^{-1}2^{-c'\varepsilon^2 N}
\end{equation}
for $\varepsilon > c\sqrt{h(G,d)}$. Combining (\ref{one}) and (\ref{two}) gives (\ref{main-prob-bound}).

We now turn to (\ref{to-show}). Observe that it is enough to prove (\ref{to-show}) for all $0 < \varepsilon' \leq \varepsilon_0$, where $\varepsilon_0 \leq 1 - a_\Lambda^+(k)$ may be any constant (perhaps depending on $H$ and $\Lambda$). Indeed, for any $\varepsilon' \geq \varepsilon_0$ we know that there is a choice of $\delta < \delta^+_k$ for which
\begin{eqnarray*}
\frac{w_{\Lambda(k,\delta)}(A^+)w_{\Lambda(k,\delta)}(B^+)}{w_\Lambda(A^+)w_\Lambda(B^+)(1+\delta)^{2(a_\Lambda^+(k)+\varepsilon')}} & \leq &  \frac{w_{\Lambda(k,\delta)}(A^+)w_{\Lambda(k,\delta)}(B^+)}{w_\Lambda(A^+)w_\Lambda(B^+)(1+\delta)^{2(a_\Lambda^+(k)+\varepsilon_0)}} \\
& \leq & 2^{-c\varepsilon_0^2}.
\end{eqnarray*}
Setting $c' = c\varepsilon_0^2$ we have $2^{-c\varepsilon_0^2} \leq 2^{-c'\varepsilon'^2}$
for $\varepsilon' \geq \varepsilon_0$ and $2^{-c\varepsilon'^2} \leq 2^{-c'\varepsilon'^2}$
for $\varepsilon' < \varepsilon_0$, so we may replace $c$ with $c'$ to obtain the result for the full range of $\varepsilon'$.
From now on we will assume that $\varepsilon' < \varepsilon_0$, for a certain $\varepsilon_0$ that will be specified later.

Setting
$$
\grg_A = \frac{\gl_k\mathbf{1}_{\{k \in A^+\}}}{2w_\Lambda(A^+)}, ~~~ \grg_B = \frac{\gl_k\mathbf{1}_{\{k \in B^+\}}}{2w_\Lambda(B^+)}
$$
(so $a_\Lambda^+(k)=\grg_A+\grg_B$) the left-hand side of (\ref{to-show}) becomes
\begin{equation} \label{split}
\frac{(w_\Lambda(A^+) + \delta \gl_k \mathbf{1}_{\{k \in A^+\}})}{(1+\delta)^{2\grg_A + \varepsilon'}w_\Lambda(A^+)} \times
\frac{(w_\Lambda(B^+) + \delta \gl_k \mathbf{1}_{\{k \in B^+\}})}{(1+\delta)^{2\grg_B + \varepsilon'}w_\Lambda(B^+)}.
\end{equation}
If either $A^+=\{k\}$ or $k \not \in A^+$ then the first term of (\ref{split}) is $(1+\delta)^{-\varepsilon'}$
so that in this case we have that for any $\delta>0$ depending only on $H$ and $\Lambda$, and any $0 < \varepsilon' \leq 1$,
\begin{equation} \label{easy-case}
\frac{(w_\Lambda(A^+) + \delta \gl_k \mathbf{1}_{\{k \in A^+\}})}{(1+\delta)^{2\grg_A + \varepsilon'}w_\Lambda(A^+)} \leq 2^{-c\varepsilon'} \leq  2^{-c\varepsilon'^2},
\end{equation}
where $c$ is a positive constant depending on $H$ and $\Lambda$ (the last inequality using $\varepsilon' \leq 1$).
If $k \in A^+$ and $|A^+|>1$ then the first term of (\ref{split}) takes the form
\begin{eqnarray}
\frac{w_\Lambda(A^+) + \delta \gl_k}{(1+\delta)^{2\grg_A + \varepsilon'}w_\Lambda(A^+)} & \leq & \frac{1 + \delta (\gl_k/w_\Lambda(A^+))}{1+\delta(2\grg_A + \varepsilon')} \nonumber \\
& = & 1 - \frac{\delta \varepsilon'}{1+\delta((\gl_k/w_\Lambda(A^+))+\varepsilon')} \nonumber \\
& \leq & 1 - \frac{\delta \varepsilon'}{3}, \label{last} 
\end{eqnarray}
with (\ref{last}) valid for sufficiently small $\varepsilon'$. Now taking $\delta = \varepsilon'$ (having chosen $\varepsilon_0$ small enough that this choice is allowed, and that (\ref{last}) holds), we get a bound of $2^{-c\varepsilon'^2}$ on the first term of (\ref{split}), where $c$ is a positive constant depending on $H$ and $\Lambda$ only. 

Repeating this analysis for the second term of (\ref{split}), we obtain (\ref{to-show}) and thus (\ref{main-prob-bound}).

\medskip

Applying (\ref{main-prob-bound}) with $\varepsilon = c\sqrt{(\log N)/N}$ (if $(\log N)/N > h(G,d)$) and $\varepsilon = c\sqrt{h(G,d)}$ (otherwise), where $c \geq c_1$ satisfies $c^2c_3 \geq 1$, we easily obtain (\ref{cor-prob-bound}), based on the observation that in both cases
$$
E_\Lambda(s(k,f)) \leq (a^+_\Lambda(k) + \varepsilon)\left(1-c_2\varepsilon^{-1}2^{-c_3\varepsilon^2N}\right) + c_2\varepsilon^{-1}2^{-c_3\varepsilon^2N}
$$
with a similar lower bound involving $a^-_\Lambda(k)$.

\section{Proof of Corollary \ref{cor-percolation}} \label{sec-perc-proof}

We assume throughout that $|V(G^n)|=N$ (a function of $n$) and that $G^n$ has fixed bipartition $\cE \cup \cO$.

We begin with the $p=\omega(1/n)$ regime. We take
$$
d= np -\sqrt{2xnp}
$$
with $x=\sqrt{f(n)}$.
The choice of $x$ is driven by the aim of making all of the terms of $h(G^n_p,d)$ be $o(1)$, with probability $1-o(1)$; this is enough for both statements of the corollary in this regime. Note that since $|\cE|=|\cO|$ we immediately have  $(|\cO|-|\cE|)/N=o(1)$.

By our choice of $x$
$$
\sqrt{2xnp} \leq \frac{np}{2}
$$
(for large enough $n$) and so $d \geq np/2$
and $1/d=o(1)$.

For a given vertex $v \in \cE$, let $d(v)$ be its degree in $G^n_p$. This is a binomial random variable with parameters $n$ and $p$, and so by standard Chernoff-type bounds (see for example \cite[Appendix A]{AlonSpencer}) we have
$$
P(d(v) < d) \leq e^{-x}.
$$
(The specific bound we are using here is
$$
P({\rm Bin}(n,p)-np < -a) < e^{-a^2/2pn}
$$
for $a>0$.)
The distribution of vertices from $\cE$ which have degree smaller than  $d$ is therefore binomial with parameters $N/2$ and $p' \leq e^{-x}$. The expected number of such vertices is at most $Ne^{-x}/2$, and by Markov's inequality the probability that there are more than $Ne^{-x+\sqrt{x}}/2$ such is at most $e^{-\sqrt{x}}$. Since $x=\omega(1)$, this is $o(1)$, and so with probability $1-o(1)$ we have
$$
\frac{|\{v \in \cE: d(v) < d\}|}{N} = o(1).
$$

It remains to consider $S:=\sum\left\{d(v)-d:v \in \cO, d(v)\geq d\right\}$. We have
\begin{eqnarray}
E(S) & = & \sum_{v \in \cO} E\left(d(v){\mathbf 1}_{\{d(v)\geq d\}}\right) - dE\left({\mathbf 1}_{\{d(v)\geq d\}}\right) \nonumber \\
& \leq & \sum_{v \in \cO} \left(\sum_{j \geq d} j{n \choose j}p^j(1-p)^{n-j}-d(1-e^{-x})\right) \nonumber \\
& \leq & N \left(np - d + de^{-x}\right) \nonumber \\
& \leq & N \left(\sqrt{2xnp} + npe^{-x}\right), \nonumber
\end{eqnarray}
and so
$$
E\left(\frac{S}{dN}\right) \leq 2\sqrt{\frac{2x}{f(n)}} + \frac{2}{e^x} ~\left(=o(1)\right)
$$
for large enough $n$ (again using $d \geq np/2$). Again by Markov's inequality, with probability $1-o(1)$ we have $S/dN=o(1)$ and so with probability $1-o(1)$ we have $h(G^n_p,d)=o(1)$, as required.

\medskip

We now deal with the $p=o(1/n)$ regime. The probability that a particular vertex is isolated in $G_p^n$ is $(1-p)^n \geq 1-2f(n)$ (for large enough $n$), so the number of non-isolated vertices in $\cE$ is a binomial random variable with parameters $N/2$ and $p' \leq 2f(n)$. By the Chernoff bound, asymptotically almost surely (with probability tending to one as $n$ tend to infinity) $\cE$ has fewer than $2f(n)N$ non-isolated vertices and so also asymptotically almost surely $G^n_p$ has fewer than $4f(n)N$ non-isolated vertices. For each $k \in V(H)$, the number of isolated vertices mapped to $k$ is a binomial random variable with parameters $m \geq N(1-4f(n))$ and $p''=\gl_k/w_\Lambda(H)$ and so (again by Chernoff bounds) asymptotically almost surely there are at least $N(1-5f(n))\gl_k/w_\Lambda(H)$ vertices of $G_p^n$ mapped to $k$. Since $\sum_{k \in V(H)} \gl_k/w_\Lambda(H) = 1$, we also have that asymptotically almost surely there are at most
$$
N\left(\frac{\gl_k}{w_\Lambda(H)} + 5f(n)\left(1 - \frac{\gl_k}{w_\Lambda(H)}\right)\right)
$$
vertices of $G_p^n$ mapped to $k$. This completes the proof of the corollary.

\section{Proof of Theorem \ref{thm-tightbounds}}
\label{sec-colourings-proof}

The graph $G_d$ will be a random $d$-regular bipartite graph on $N=c^{d/\log d}$ vertices (where $c>1$ will depend on the particular $H$ and $\Lambda$ under consideration). A standard method of constructing such a graph is as follows. We begin with a set of size $Nd$ consisting of $Nd/2$ type I vertices $\{u_{ij}:1 \leq i \leq N/2, 1 \leq j \leq d\}$ and $Nd/2$ type II vertices $\{v_{ij}:1 \leq i \leq N/2, 1 \leq j \leq d\}$.
We then choose a uniformly random perfect matching from the type I vertices to the type II vertices, and turn this into a $d$-regular bipartite multigraph on $N$ vertices with bipartition classes $\cE=\{u_1, \ldots, u_{N/2}\}, \cO=\{v_1, \ldots, v_{N/2}\}$ by, for each $i=1, \ldots, N/2$, identifying $u_{i,1}, \ldots, u_{i,d}$ with $u_i$ and $v_{i,1}, \ldots, v_{i,d}$ with $v_i$. Finally, we condition on the result being a simple graph. This process generates a $d$-regular bipartite graph on $N$ vertices with bipartition classes $\cE$, $\cO$, uniformly (see for example \cite{Wormold}).

O'Neil \cite{Oneill} showed that the probability that the multigraph produced by this process is simple is (for large enough $d$) at least $e^{-d^2/3}$. It follows that if we establish that the multigraph produced (before conditioning on being simple) has a certain property with probability at least $1-e^{-d^2}$ (say), then there is a simple $d$-regular graph with that property.

We want to establish that for large enough $d$ the multigraph has a number of desirable expansion properties. First, we want to show that for each $C\log d  \leq j \leq 3\log dN/d$
(for some constant $C>0$, depending on $c$),
every subset of $\cE$ of size $j$ and every subset of $\cO$ of size $j$ has at least $\alpha j$ distinct neighbours where $\alpha = d/(C\log d)$. For a particular such $j$, the probability that the graph fails to have this property is (by a union bound) at most
\begin{eqnarray*}
2{N/2 \choose j}{N/2 \choose \alpha j} \frac{(\alpha jd)_{jd}}{(Nd/2)_{jd}} & \leq & \left(\frac{eN}{2 \alpha j}\right)^{2 \alpha j}\left(\frac{2\alpha j}{N}\right)^{jd} \\
& = & e^\frac{2jd}{C\log d} \left(\frac{2jd}{CN \log d}\right)^{jd - \frac{2jd}{C\log d}} \\
& \leq & e^\frac{2jd}{C\log d} \left(\frac{2jd}{CN \log d}\right)^{jd/2}
\end{eqnarray*}
(for large enough $d$, depending on $C$) with the first inequality using ${n \choose r} \leq (en/r)^r$. For $j \geq d\log d$ we bound $2jd/(CN\log d) \leq  1/2$ (valid for $C \geq 12$) so that for large enough $d$ (depending on $C$)
$$
e^\frac{2jd}{C\log d} \left(\frac{2jd}{CN \log d}\right)^{jd/2} \leq 1.4^{-jd} \leq e^{-2d^2}.
$$
For $j \leq d \log d$ we instead bound $(2dj)/(C N \log d) \leq d^2/N$ (valid for $C \geq 2$). We now have
$$
e^\frac{2jd}{C\log d} \left(\frac{2jd}{CN \log d}\right)^{jd/2} \leq \exp\left\{2jd\log d - \frac{jd^2 \log c}{2 \log d}\right\} \leq \exp\left\{\frac{-jd^2 \log c}{3 \log d}\right\}
$$
(again for large $d$, recalling $N=c^{d/\log d}$), which is at most $e^{-2d^2}$ for $j \geq C\log d$ for suitable $C$ depending on $c$. Since there are at most $N=c^{d\log d}$ choices for $j$, the probability that the graph fails to have the desired property for some $j$ is at most $e^{-d^2}$. If the process results in a simple graph, then we trivially get the same expansion for subsets of $\cE$ or $\cO$ of size at most $C\log d$, since for $1 \leq j \leq C\log d$ there is a trivial lower bound of $d$ on the neighbourhood size of a set of size $j$, and we have $d \geq jd/(C\log d)$ for $j$ in this range.

Next we establish that the graph has the property that for every subset $A$ of $\cE$ of size $3N\log d/d$ and every subset $B$ of $\cO$ of size $3N\log d/d$, there is an edge joining a vertex of $A$ to a vertex of $B$. By a union bound, the probability that the multigraph fails to have the property is at most
$$
{N/2 \choose \beta N}^2 \frac{\left(Nd/2-\beta Nd\right)_{\beta Nd}}{\left(Nd/2\right)_{\beta Nd}} \leq \exp\left\{2\beta N\log (e/(2\beta))-2\beta^2 d N\right\}
$$
where $\beta = 3\log d/d$. With $N=c^{d/\log d}$, this is at most $e^{-d^2}$ for large enough $d$ (depending on $c$).
We have shown the following.
\begin{lemma} \label{lem-exp_of_rand_reg_graph}
Fix $c>1$. There are $d_0\geq 1$ and positive $C$, both depending on $c$, such that for all $d\geq d_0$ there is a $d$-regular,
bipartite graph $G_d$ on $N=c^{d/\log d}$ vertices with bipartition classes $\cE$ and $\cO$ satisfying the following:
\begin{enumerate}
\item Every subset of $\cE$ or $\cO$ of size $j$, with $1 \leq j \leq 3N\log d/d$, has at least $jd/(C \log d)$ neighbours.
\item Every pair of subsets each of size $3N\log d/d$, one from $\cE$ and one from $\cO$, have an edge between them.
\end{enumerate}
\end{lemma}
We now fix such a $G_d$ and study $Z_\Lambda(G_d, H)$.
Given $f \in {\rm Hom}(G_d,H)$ set
$$
\cE(f) = \{k \in V(H): |f^{-1}(k) \cap \cE| \geq 3N\log d/d\} 
$$
and
$$
\cO(f) = \{k \in V(H): |f^{-1}(k) \cap \cO| \geq 3N\log d/d\}. 
$$
Clearly both $\cE(f)$ and $\cO(f)$ are non-empty, and by Lemma \ref{lem-exp_of_rand_reg_graph}, we have $\cE(f) \sim \cO(f)$ (that is, everything in $\cE(f)$ is adjacent to everything in $\cO(f)$). So we can partition ${\rm Hom}(G_d,H)$ into classes indexed by pairs $(A,B)$ with $A \sim B$. Write ${\mathcal C}(A,B)$ for the class corresponding to $(A,B)$.
We want to establish that for $(A,B) \in \cM_\Lambda(H)$ we have
\begin{equation} \label{ref-extreme-partition}
\sum_{f \in {\mathcal C}(A,B)} w_\Lambda(f) = (1+o(1))\eta_\Lambda(H)^{N/2}
\end{equation}
while for all other $(A,B)$ we have
\begin{equation} \label{ref-less-extreme-partition}
\sum_{f \in {\mathcal C}(A,B)} w_\Lambda(f) = o\left(\eta_\Lambda(H)^{N/2}\right),
\end{equation}
where all asymptotic terms are (unless stated otherwise) as $d \to \infty$. 
From this we see that
$$
Z_\Lambda(G_d,H) = |\cM_\Lambda(H)|(1+o(1))\eta_\Lambda(H)^{N/2},
$$
and that all but a vanishing proportion of $Z_\Lambda(G_d,H)$ comes from pure-$(A,B)$ colourings (with $(A,B)\in \cM_\Lambda(H)$) in which $\cE$ is mapped to $A$ and $\cO$ to $B$, with each such $(A,B)$ contributing equally to $Z_\Lambda(G_d,H)$; this is enough to give the first part of Theorem \ref{thm-tightbounds}. Indeed, fix $(A,B) \in \cM_\Lambda(H)$. A proportion $(1+o(1))/|\cM_\Lambda|$ of $Z_\Lambda(G_d,H)$ is obtained by independently colouring $\cE$ from $A$ and $\cO$ from $B$ according to the given weights. Fix $k \in A$. We claim that with very high probability, a proportion very close to $\gl_k/w_\Lambda(A)$ of $\cE$ gets mapped to $k$. Set $p=\gl_k/w_\Lambda(A)$ and $m=N/2$. The number $U_k$ of vertices of $\cE$ mapped to $k$ is a binomial random variable with parameters $m$ and $p$. So by Tchebychev's inequality,
$$
{\rm Pr}\left(|U_k-p m| \geq \log m \sqrt{m p (1-p)}\right) \leq \frac{1}{\log^2 m}.
$$
This shows that the proportion of vertices mapped to $k$ in a pure-$(A, B)$ colouring is very close to
$$
\frac{\gl_k \mathbf{1}_{\{k \in A\}}}{2w_\Lambda(A)} + \frac{\gl_k \mathbf{1}_{\{k \in B\}}}{2w_\Lambda(B)}
$$
with high probability. Applying this with $(A,B)=(A^+,B^+)$ and $(A,B)=(A^-,B^-)$, the first part of Theorem \ref{thm-tightbounds} follows.

The lower bound in (\ref{ref-extreme-partition}) is obtained by considering pure-$(A,B)$ colourings with $\cE$ mapped to $A$ and $\cO$ to $B$. To establish (\ref{ref-less-extreme-partition}) and the upper bound in (\ref{ref-extreme-partition}), fix $0 \leq j \leq 3N\log d/d$ let $q = |V(H)|$, and assume that $d$ is large. 
We consider the contribution to $\sum_{f \in {\mathcal C}(A,B)} w_\Lambda(f)$ from those $f \in {\mathcal C}(A,B)$ in which, for each $k \not \in A \cup B$, we have at most $j$ vertices mapped to $k$, and we have at least one $k' \not \in A \cup B$ whose preimage has size $j$. 
To bound the contribution from these $f$, we first bound the number of ways of locating the vertices that are mapped to $k$ for each $k \notin A \cup B$ by 
$\left(\sum_{i \leq j} {N \choose i}\right)^q$.
The contribution to the sum of the weights from these exceptional vertices is at most
$w_\Lambda(H)^{qj}$.
For the contribution from the remaining vertices, we deal separately with the cases $(A,B) \in \cM_\Lambda(H)$ and $(A,B) \not \in \cM_\Lambda(H)$. For $(A,B) \not \in \cM_\Lambda(H)$, we simply upper bound the contribution by $\left(w_\Lambda(A)w_\Lambda(B)\right)^{N/2}$, leading to
\begin{eqnarray*}
\sum_{f \in {\mathcal C}(A,B)} w_\Lambda(f) & \leq & \left(w_\Lambda(A)w_\Lambda(B)\right)^\frac{N}{2}   \left(\sum_{i \leq j} {N \choose i}\right)^q (w_\Lambda(H))^{qj} \\ 
& = & o\left(\eta_\Lambda(H)^{N/2}\right),
\end{eqnarray*}
as required.
For $(A,B) \in \cM_\Lambda(H)$, consider a $k'$ that has preimage size $j$. We claim that there are at least $jd/(2C\log d)$ vertices which, in the specification of $f$, need to be mapped to $A \cup B$ and which are adjacent to at least one of the $j$ vertices mapped to $k'$. Indeed, by Lemma \ref{lem-exp_of_rand_reg_graph}, the neighbourhood size of the $j$ vertices mapped to $k'$ is at least $jd/(C\log d)$, and at most $qj$ vertices have been mapped to vertices from outside $A \cup B$, so there are at least $jd/(C\log d)-qj > jd/(2C\log d)$ vertices that are adjacent to a vertex mapped to $k'$ and need to be mapped to vertices from $A \cup B$. Since $k'$ cannot be adjacent to everything in $A$, nor can it be adjacent to everything in $B$ (else we would not have $(A, B) \in \cM_\Lambda(H)$), our choice on these at least $jd/(2C\log d)$ vertices is restricted to a proper subset of $A \cup B$; the contribution we get from the remaining vertices (those mapped to $A \cup B$) is therefore at most
$$
\frac{\left(w_\Lambda(A)w_\Lambda(B)\right)^{\frac{N}{2}}}{(1+\varepsilon)^{\frac{jd}{2C\log d}}}
$$
where $\varepsilon > 0$ (depending on $H$ and $\Lambda$) can be chosen uniformly for all $A$, $B$. Combining these observations we get that
$$
\sum_{f \in {\mathcal C}(A,B)} w_\Lambda(f) \leq \eta_\Lambda(H)^\frac{N}{2}  \frac{ \left(\sum_{i \leq j} {N \choose i}\right)^q (w_\Lambda(H))^{qj}}{(1+\varepsilon)^{\frac{jd}{2C\log d}}}. 
$$
If $j=0$, the right-hand side above is $(w_\Lambda(A)w_\Lambda(B))^{N/2}$. For $j > 0$ it can be bounded above by
$$
\left(\frac{1}{(1+\varepsilon')^{\frac{d}{\log d}}}\right)^j
$$
for some $\varepsilon' > 0$ (depending on $H$ and $\Lambda$) for all $j$ in the range $1 \leq j \leq 3N\log d/d$, as long as $c$ is sufficiently small (recall $N=c^{d/\log d}$). Summing over $j$ gives the upper bound in (\ref{ref-extreme-partition}).

\medskip

We now turn to the second part of Theorem \ref{thm-tightbounds}. We take $G'_d$ to be the disjoint union of $m$ copies of $K_{d,d}$ where $m =m(d) = \omega(1)$. Fix $k \in V(H)$. Let $X$ be the number of vertices mapped to $k$ in a $p_\Lambda$-chosen $H$-colouring of $G'_d$, and $X_i$ the number mapped to $k$ in the $i$th copy of $K_{d,d}$. Define $a_\Lambda(k)$ by $E(X_i)=2d a_\Lambda(k)$, and note that ${\rm Var}(X_i) \leq d^2$. Since $X=\sum_{i=1}^m X_i$ we have $E(X)=2dma_\Lambda(k)$ and ${\rm Var}(X) \leq md^2$. By Tchebychev's inequality,
$$
P(|X-2dma_\Lambda(k)| > 2dm\varepsilon) = P(|X/2dm-a_\Lambda(k)|> \varepsilon) \leq 1/4m\varepsilon^2.
$$
So choosing $\varepsilon = o(1)$ with $m\varepsilon^2 = \omega(1)$ (for example, $\varepsilon = 1/m^{1/3}$), the probability that the proportion of vertices mapped to $k$ in a $p_\Lambda$-chosen $H$-colouring of $G'_d$ differs from $a_\Lambda(k)$ by more than $o(1)$ is at most $o(1)$. The claimed bound on $s(k,f)$ follows, as does the estimate of $\bar{p}_\Lambda(k)$.

\end{document}